\documentclass[11pt]{article}
\PassOptionsToPackage{numbers, compress}{natbib}
\usepackage[preprint]{neurips_2021}
\usepackage{amsmath}
\usepackage{amssymb}
\usepackage[utf8]{inputenc} 
\usepackage[T1]{fontenc}    
\usepackage{hyperref}       
\usepackage{url}            
\usepackage{booktabs}       
\usepackage{amsfonts}       
\usepackage{nicefrac}       
\usepackage{microtype}      
\usepackage{xcolor}         
\setcitestyle{authoryear,round}
\makeatletter
\usepackage{amsthm}\usepackage{amsfonts}
\usepackage{bbm}\usepackage{epsfig}
\newtheorem{theorem}{Theorem}
\newtheorem{definition}{Definition}
\newtheorem{corollary}{Corollary}
\newtheorem{lemma}{Lemma}

\usepackage{listings}

\newcommand{\Rmnum}[1]{\expandafter\@slowromancap\romannumeral #1@}

\usepackage{authblk}

\begin{document}
\title{Asymptotic Inference for Infinitely Imbalanced Logistic Regression}
\author[1]{Dorian Goldman}
\author[2]{Bo Zhang}
\affil[1]{Lyft, dgoldman@lyft.com}
\affil[2]{bozhang@live.com}

\maketitle 
\begin{abstract}
    In this paper we extend the work of \cite{owen} by deriving a second order expansion for the slope parameter in logistic regression, when the size of the majority class is unbounded and the minority class is finite. More precisely, we demonstrate that the second order term converges to a normal distribution and explicitly compute its variance, which surprisingly once again depends only on the mean of the minority class points and not their arrangement under mild regularity assumptions. In the case that the majority class is normally distributed, we illustrate that the variance of the the limiting slope depends exponentially on the z-score of the average of the minority class's points with respect to the majority class's distribution. We confirm our results by Monte Carlo simulations.
\end{abstract}

\section{Introduction}

Logistic regression is one of the most widely used machine learning algorithms for binary classification problems. In many binary classification problems, one of the classes is very rare. Common examples include fraud detection (\cite{bolton2002statistical}), ad conversion (\cite{lee2012estimating}) and drug interactions (\cite{zhu2006lago}). This paper studies the infinitely imbalanced regime originally formulated and analyzed by \citet{owen} and derives a closed-form estimator for the logistic regression coefficient's variance, and hence its confidence interval. Our result is obtained through a series expansion approach and can be viewed as a refinement to the first-order limit result in \cite{owen}.

Specifically, we consider a data set of $(x,y)$ pairs with $x\in \mathbb{R}^d$ and $y\in \{0,1\}$, consisting of $n$ observations with $y=1$ and $N$ observations with $y=0$, where $n$ is finite and $N \to +\infty$. Denote by $\beta_N\in\mathbb{R}^d$ the vector of regression coefficients resulting from applying logistic regression to the data set. Our main result is that
\begin{equation}\label{eq:intro_main_result}
\sqrt{N}(\beta_N - \beta_*) \overset{\mathcal{D}}{\rightarrow} \mathcal{N}\left(0, \Sigma \right) , \textrm{ as } N \to +\infty,
\end{equation}
where $\overset{\mathcal{D}}{\rightarrow}$ means convergence in distribution, $\beta_*\in\mathbb{R}^d$ satisfies
\begin{equation}\label{eq:owen_main00}
\frac{\int e^{\beta_* \cdot x^T}x dF_0(x)}{\int e^{\beta_* \cdot x^T} dF_0(x)} = \bar x,
\end{equation}
we denote by superscript $^T$ the transpose of a vector, $F_0$ is the cdf (cumulative distribution function) for the majority class, $X \lvert Y = 0$, $\bar{x}$ is the average of the sample $x_i$ values for which $y=1$, and the asymptotic covariance $\Sigma$ is given by
\begin{align}
    \Sigma = \Sigma_{\bar x, \mu}&:= H(\bar x, F_0)^{-1} V(\bar x,F_0) \left(H(\bar x, F_0)^{-1}\right)^T,
\end{align}
with
\begin{align}
  H(\bar x, F_0) &:= \mathbb{E}_{F_0}\left[e^{\beta_* X^T}(X - \bar X)^T (X - \bar X)\right],\\
    V(\bar x, F_0) &:= \mathbb{E}_{F_0}\left[e^{2\beta_* X^T}(X - \bar X)^T (X - \bar X)\right].
\end{align}
 In the one-dimensional case, the above expression for $\Sigma_{\bar x, \mu}$ reduces to a scalar,
\begin{equation}\label{eq:result_1d}
     \Sigma_{\bar x, \mu} = \sigma_{\bar x, \mu}^2 := \frac{ \int e^{2\beta_*x} (x - \bar x)^2 dF_0(x)}{\left(\int e^{\beta_*x} (x - \bar x)^2 dF_0(x)\right)^2}.
\end{equation}
Our result suggests $(\beta_*- Az_{\theta/2}/ \sqrt{N}, \beta_*+ Az_{\theta/2}/ \sqrt{N})$ as the $100(1-\theta)$ percent confidence interval estimator of $\beta_N$, 
where $A$ is the Choleski decomposition matrix of $\Sigma$ and $z_{\theta/2}$ is the $(1-\theta/2)$ percentile of the standard normal distribution.
Also, the estimator is asymptotically exact as $N\rightarrow\infty$.

We do not assume that the data generating process is truly specified by the Logistic Regression model as a Generalized Linear Model, i.e., that the response variables are independent Bernoulli random variables with each of their log-odds being the same linear combination of the independent variables. The only related inference result under model misspecification in the literature is that of the \textit{sandwich estimator} (e.g., see Section 8.3.2 of \cite{agresti2015foundations}). Specifically, an alternative expression for the variance in the forementioned one-dimensional case is
\begin{equation}
   \sigma_{\bar x,\mu}^2= \frac{\textrm{Var}_{F_0}(\partial_{\beta} g(\cdot ; \beta=\beta_*))}{\mathbb{E}_{F_0}[\partial_{\beta}^2 g(\cdot ; \beta=\beta_*)]^2}, 
\end{equation}
where $g(x ; \beta) = e^{\beta(x-\bar x)}$. This expression has the same form as the sandwich estimator of the maximum likelihood estimator's variance under misspecification (see Section 16.2 of \cite{ZF2016}).




Our result can be viewed as a refinement of the limit result in \cite{owen}. Specifically, \cite{owen} shows that \begin{equation}\label{eq:owen_main}
    \lim_{N \to +\infty}\frac{\int e^{\beta_N \cdot x^T}x dF_0(x)}{\int e^{\beta_N \cdot x^T} dF_0(x)} = \frac{\int e^{\beta_* \cdot x^T}x dF_0(x)}{\int e^{\beta_* \cdot x^T} dF_0(x)} = \bar x.
\end{equation}
That is, in the infinitely imbalanced limit, the vector of regression coefficients is the tilting parameter such that the first moment of the majority class feature's exponentially tilted distribution matches with the sample mean of the minority class's feature vector. The remarkable aspect of the result is that the regression coefficients in the limit depend only on the average of the points $x$ where $y=1,$ and not on how these points are distributed. The limiting slope $\beta_*$ is unique since it is the solution to a strictly convex problem when $H$ is invertible. 



 
Our result enables an analytical assessment of the accuracy of approximating $\beta_N$ by $\beta_*$. Consider the example of $F_0$ being the cdf of $\mathcal{N}(\mu, \sigma^2)$. In this case, the limiting variance
\begin{equation}
    \sigma_{\bar x, \mu}^2= e^{\frac{(\bar x - \mu)^2}{\sigma^2}}\left[\frac{(\bar x - \mu)^2 + \sigma^2}{\sigma^4}\right],
\end{equation}
which depends exponentially on the number of standard deviations of $\bar x$ to $\mu$, or the z-score of $\bar x$ with respect to $\mathcal{N}(\mu, \sigma^2)$. As a result, in order for $\beta_*$ to accurately estimate $\beta_N$, an exceedingly large number of data points can be required if $\frac{|\bar x - \mu|}{\sigma}$ is large. Specifically, since \begin{equation}
    \beta_N\overset{\mathcal{D}}{=} \mathcal{N}(\beta_*, \sigma^2_{\bar x, \mu}N^{-1})+o(N^{-1/2}),
\end{equation}
it holds approximately that
\begin{equation}
    \beta_N\in [\beta_*-\sigma_{\bar x, \mu}N^{-1/2}z_{\theta/2}, ~\beta_*+\sigma_{\bar x, \mu}N^{-1/2}z_{\theta/2}], \mbox{~with probability $1-\theta$.}
\end{equation}
In order to achieve $|\beta_N - \beta_*| < \epsilon$, we need 
\begin{equation}
    N \sim \frac{e^{2(\bar x - \mu)^2/\sigma^2}}{\epsilon^2}.
\end{equation} 
For example, if $\bar{x}$ is more than three-$\sigma$ away from $\mu$, then for a desired accuracy of $10^{-1}$ one would need roughly $10^8$ samples. This is not feasible in many practical cases.

The remainder of this paper is organized as follows. Section 2 collects needed notation and assumptions. In section 3 we prove our main results and explore the special cases when $F_0$ is the normal distribution.  In Section 4 we provide numerical simulations which confirm the expansion and illustrate through examples where it degrades.

\section{Notation and Assumptions}

\begin{itemize}
\item We denote sampled vectors in lower case $z \in \mathbb{R}^d$ where $d$ is the dimension and their corresponding random variable by its capitalization $Z$. We assume by default that samples are row vectors, ie. $x = (x_1,x_2,\cdots,x_d)$ with $x_i \in \mathbb{R}^n$ for $i = 1,2,\cdots,n$  and hence the dot products between two vectors $x$ and $z$ is denoted $x z^T$. 
\item  Let $(\Omega, \mathcal{F}_0,P)$ be the probability space on infinite majority samples $X \lvert Y=0$. When we say \emph{almost surely}, this refers to the probability measure induced on $\mathcal{F}_0(\Omega)$, the sigma algebra of $P-$measurable sets in $\Omega$. See \cite{billingsley1999convergence} for a detailed construction.

\item We assume that $ X \lvert Y=0$ has pdf $F_0$ which is absolutely continuous with respect to Lebesgue measure on $\mathbb{R}^d$ (ie. $dP = F_0dx$) and that $X \lvert Y=1$ is a finite collection of points with empirical mean $\bar x$. We write $x_i \sim F_0$ to denote that $x_i$ is sampled from $F_0$. 

\item There are $N$ observations with $Y=0$ denoted  $\{x_i\}_{i=1}^N$ and $n$ observations with $Y=1$ denoted $\{x_j\}_{j=1}^n$.

\item The expectation of the random variable $Z$ with respect to a probability distribution $G$ is denoted $\mathbb{E}_{G} [Z]$. 

\item We assume that $F_0$ has finite first and second tilted moments for all $\beta \in \mathbb{R}^d$:
\begin{equation}\label{eq:finite_moments}
\int e^{\beta \cdot x^T }(1 + \|x\|^k) dF_0(x) < +\infty \textrm{ for } 1 \leq k \leq 2.
\end{equation}
This is a stronger condition than \cite{owen} assumes when obtaining the first order result.
\item We also assume the range of the finite points $\{x_j\}_{j=1}^n$ is bounded, i.e., \begin{equation}\label{eq:xjbound}\sup_{N>0} \max_{1 \leq j \leq n} \|x_j\| < C < +\infty.\end{equation}
We believe the above assumption is not necessary, and removing it  will simply add an additional term into our asymptotic variance. However we make this assumption to simplify the arguments. 
\end{itemize}

We center our log-loss around $\bar x$ as in \cite{owen}: 
\begin{equation}\label{eq:logloss_empirical2}
    \mathcal{L}_N(\alpha, \beta) = n\alpha -  \sum_{j=1}^n\log(1 + e^{\alpha+\beta(x_j-\bar x}) - \sum_{i=1}^N \log (1 + e^{\alpha + \beta (x_i - \bar x)}).
\end{equation}

It is well known that when the two classes are perfectly separated by a hyperplane, no solution exists to the logistic regression problem. Indeed, a degree of overlapping is required to ensure existence of finite solutions to the problem. This problem was fully explored in \cite{Silvapulle1981}, where he completely characterizes the existence and uniqueness criteria for the logistic regression problem. \cite{owen} makes a slightly stronger assumption than \cite{Silvapulle1981}, specifically, that $F_0$ \emph{surrounds} the point $\bar x$ in the following sense:

\begin{definition}\label{def:surrounds}
We say that $F_0$ surrounds $\bar x$ if there exists some $\epsilon > 0$ and $\delta > 0$ such that for all $\omega$
\begin{equation}
   \mathop{\inf_{\omega}}_{\omega^T \omega=1}  \int_{(x-\bar x) \cdot \omega > \epsilon } dF(x) > \delta.
\end{equation}
\end{definition}

The above definition ensures that if we consider any possible hyperplane intersecting $\bar x$, then $F_0$ will assign positive mass to a ball containing $\bar x$ intersected with the half plane $(x - \bar x) \cdot w > \epsilon$. \\

\section{Main Results} 

Our main result is as follows.

\begin{theorem} \label{theorem:main}
Assume that $F_0$ satisfies \eqref{eq:finite_moments} and surrounds $\bar x$ in the sense of Definition \ref{def:surrounds} and that \eqref{eq:xjbound} holds. Then there exists a unique solution, denoted $\beta_*$,  to \eqref{eq:owen_main} and a unique maximizer $(\beta_N,\alpha_N)$ of \eqref{eq:logloss_empirical2}. Additionally, assume that \begin{equation}\label{eq:det_lbound}\det(\mathbb{E}_{F_0}[e^{\beta_*X^T}(X-\bar x)^T (X - \bar x)])>0.\end{equation}
 Then it follows that
\begin{equation}
    \lim_{N \to +\infty} \sqrt{N}(\beta_N - \beta_*)  \overset{\mathcal{D}}{=} \mathcal{N}\left(0, \Sigma\right)\textrm{ as } N \to +\infty,
\end{equation}
where 
\begin{align}
    \Sigma &:= H(\bar x, F_0)^{-1} V(\bar x,F_0) \left(H(\bar x, F_0)^{-1}\right)^T\\ 
    H(\bar x, \mu) &:= \mathbb{E}_{F_0}\left[e^{\beta_* X^T}(X - \bar X)^T (X - \bar X)\right]\\
    V(\bar x, \mu) &:= \mathbb{E}_{F_0}\left[e^{2\beta_* X^T}(X - \bar X)^T (X - \bar X)\right].
\end{align}

\end{theorem}

The following corollary provides an explicit estimate for the variance in the case that $F_0$ is a 1D Normal distribution. Recall that when $F_0$ is the cdf of $\mathcal{N}(\mu,\sigma^2)$, we have $\beta_*=\frac{\bar x - \mu}{\sigma^2}$ from \cite{owen}. Our main result in this case allows us to conclude the following corollary:\\

\begin{corollary}\label{cor:1d_normal}
Assume in addition to the assumptions of Theorem \ref{theorem:main} that $F_0$ is a normal distribution with mean $\mu$ and variance $\sigma^2$. Then  $\Sigma_{\bar x \mu}$ can be represented by a scalar $\sigma_{\bar x, \mu}^2$ defined as
\begin{equation}
     \sigma_{\bar x \mu}^2 := e^{\frac{(\bar x - \mu)^2}{\sigma^2}}\left[\frac{(\bar x - \mu)^2 + \sigma^2}{\sigma^4}\right].
\end{equation}
\end{corollary}
\noindent \textbf{Proof:} Use equation \ref{eq:result_1d} and then see derivation after Corollary \ref{cor:variance_ratio} in the Appendix. \\

There are two interesting observations of Corollary \ref{cor:1d_normal}. One is that the variance, and hence uncertainty, depends only on the distance of $\bar x$ to $\mu$. The second,  is that the variance of the estimator grows exponentially with the number of standard deviations $\bar x$ is from the mean of $F_0$. This implies that the first order limiting solution deteriorates very quickly when $\bar x$ is far from the center of mass of $F_0$. Our simulations in Section \ref{section:simulations} confirm this implication.

The primary application of Theorem \ref{theorem:main} is that it provides more accurate inference of the slope parameter in the highly imbalanced case, without needing to assume that the data generation process is consistent with the assumptions of Logistic Regression (See Section 4).

In this section we will prove Theorem 1. The idea of the proof is quite simple - perform a Taylor expansion about the limiting solution as one does with MLE estimators. However the nuance here is carefully estimating the decaying effect of $e^{\alpha_N}$ to extract the limiting behavior of the functional.

We begin by repeating a result from \cite{owen}, which gives us the desired decay of the intercept parameter $\alpha$. \\

\begin{lemma}\label{lemma:alpha_conv}
Let $F_0$ surround $\bar x$ and $(\beta_N,\alpha_N)$ be maximizers of \eqref{eq:logloss_empirical2}. Then $e^{\alpha_N} \leq 2\frac{n}{N\delta}$  and $\limsup_{N \to +\infty} \|\beta_N\| <+\infty$
\end{lemma}

We now wish to use Lemma \ref{lemma:alpha_conv} to control the derivatives of the log likelihood $\mathcal{L}_N$. This will be used in the proof of Theorem 1 when passing to the limit as $N \to +\infty$.\\

\begin{lemma}\label{lemma:sum_approx}
Assume the conditions of Theorem 1. Then for a.s every sequence $\{x_i\}_{i=1}^N \sim F_0$ we have
\begin{align}\label{eq:der1_difference}
 \limsup_{N \to +\infty} e^{-\alpha_N}\left|\frac{1}{N} \sum_{i=1}^N \frac{e^{\beta_*(x_i - \bar x)}}{(1+ e^{\alpha_N+\beta_* (x_i-\bar x)})}(x_i^k-\bar x^k)\right. -  &\left. e^{\beta_*^T(x_i - \bar x)}(x_i^k-\bar x^k)\right| \leq +\infty
\end{align}

\begin{align}\label{eq:der2_difference}
\nonumber \limsup_{N \to +\infty} e^{-\alpha_N}\left|\frac{1}{N} \sum_{i=1}^N \frac{e^{\xi_N^T(x_i - \bar x)}}{(1+ e^{\alpha_N+\xi_N^T (x_i-\bar x)})^2}(x_i^k-\bar x^k) (x_i^j - \bar x^j)\right. \nonumber -  &\left.e^{\xi_N^T(x_i - \bar x)}(x_i^k-\bar x^k)(x_i^j-\bar x^j)\right|\\
&<+\infty 
\end{align}

\end{lemma}
\noindent \textbf{Proof:} 
We prove for the first inequality, the second follows similarly. Subtracting the two sums in the left side of \eqref{eq:der1_difference} we have
\begin{align}
   \nonumber \frac{1}{N}\sum_{i=1}^N \frac{e^{\beta_*^T(x_i - \bar x)}}{1+ e^{\alpha_N+\beta_* (x_i-\bar x)}}(x_i-\bar x) &- \frac{1}{N}\sum_{i=1}^N e^{\beta_*^T(x_i - \bar x)}(x_i-\bar x)\\ &=  \frac{e^{\alpha_N}}{N}\sum_{i=1}^N \frac{e^{2\beta_*^T(x_i - \bar x)}}{1+ e^{\alpha_N+\beta_*^T (x_i-\bar x)}}(x_i-\bar x).\label{eq:sum_approx}
\end{align}
It is easily seen that 
\begin{equation}\label{eq:bound1}
    \left|\frac{e^{2\beta_*^T(x_i - \bar x)}}{1+ e^{\alpha_N+\beta_*^T (x_i-\bar x)}}(x_i-\bar x) \right| \leq e^{2\beta_*^T(x_i - \bar x)}\|x_i-\bar x\|\leq e^{2\beta_*^T(x_i - \bar x)}(\|x_i\|+\|\bar x\|).
\end{equation}

By \eqref{eq:finite_moments}, the L.L.N implies that almost surely

\begin{equation}\label{eq:der1_lln}
    \frac{1}{N}\sum_{i=1}^N e^{2\beta_*^T(x_i - \bar x)}(\|x_i\|+\|\bar x\|) \to \mathbb{E}_{F_0}[e^{2\beta_*^T(x - \bar x)}(\|x\|+\|\bar x\|)].
\end{equation}

Now applying \eqref{eq:der1_lln} to \eqref{eq:bound1}, we have for sufficiently large $N$, that a.s 

\begin{align}
    \left|\frac{1}{N}\sum_{i=1}^N \frac{e^{2\beta^T(x_i - \bar x)}}{1+ e^{\alpha+\beta^T (x_i-\bar x)}}(x_i-\bar x)\right| &\leq \frac{1}{N}\sum_{i=1}^N e^{2\beta^T (x_i-\bar x)}|x_i| + |\bar x| \frac{1}{N}\sum_{i=1}^N e^{2\beta^T(x_i-\bar x)}\\
    &\leq 2\mathbb{E}_{F_0}[e^{2\beta^T(x - \bar x)}(|x|+|\bar x|)].
\end{align}
Taking limsups on both sides of the above equation and using \eqref{eq:der1_lln} yields the result. The proof of \eqref{eq:der2_difference} follows similarly.
$\Box$

Now we show convergence of the Hessian:





\begin{lemma}\label{lemma:der_bounds}
Assume once again the setup of Theorem 1 and assume $\xi_N \to \beta_*$ a.s as $N \to +\infty$. Then a.s we have 

\begin{align}\label{eq:2ndder_conv}
 \nonumber \lim_{N \to +\infty} \frac{1}{N} \sum_{i=1}^N \frac{e^{\xi_N^T(x_i - \bar x)}}{(1+ e^{\alpha_N+\xi_N^T (x_i-\bar x)})^2}(x_i^k-\bar x^k) (x_i^j - \bar x^j) \to  \mathbb{E}_{F_0} [e^{(\beta_*)^T(x_i - \bar x)}(x^k-\bar x^k)(x^j-\bar x^j)] .\nonumber
\end{align}
\end{lemma}
\noindent \textbf{Proof:} 
Define 
\begin{equation}
    G(M,N) := \frac{1}{N}\sum_{i=1}^N e^{\xi_M^T(x_i - \bar x)}(x_i^k-\bar x^k)(x_i^j - \bar x^j).
\end{equation}

We wish to show that the limit $\lim_{(M,N) \to +\infty} G(M,N)$ exists. As such we check the criteria of the Moore-Osgood theorem \cite{moore-osgood}:
\begin{itemize}
    \item 1. For each $N$, the limit $\lim_{M \to + \infty}G(M,N)$ exists pointwise. Indeed:
    \begin{align}\lim_{M \to + \infty}G(M,N) &= \lim_{M \to + \infty}\frac{1}{N}\sum_{i=1}^N e^{\xi_M^T(x_i - \bar x)}(x_i^k-\bar x^k)(x_i^j - \bar x^j)\\ &=  \frac{1}{N}\sum_{i=1}^N e^{\beta_*^T(x_i - \bar x)}(x_i^k-\bar x^k)(x_i^j - \bar x^j,)\end{align}
    which follows from \eqref{eq:finite_moments}, and the L.L.N. 
    \item 2. The limit $\lim_{N \to +\infty} G(M,N)$ converges uniformly with respect to $M$. This requires a bit more work which we show below. 
\end{itemize}

 We would like to apply the the Uniform Law of Large Numbers to $G(\cdot, N)$ and this requires showing 
\begin{itemize}
\item 1. $\{\xi_M\}_M$ is compact.
\item 2. $x \mapsto h(x,\xi_M) := e^{\xi_M^T(x_i - \bar x)}(x_i^k-\bar x^k)(x_i^j - \bar x^j)$ is continuous for each $\xi_M$ and
\item 3. There exists a dominating function $d(x)$ such that 
\begin{equation}
     h(x,\xi_M) \leq d(x) \in L_{F_0}^1(\Omega).
\end{equation}
\end{itemize}
Requirement 1 is easily satisfied since we have  $\xi_M \to \beta^*$ as $M \to +\infty$ almost surely. Consequently there exists a compact subset $K \subset \mathbb{R}^d$ with $\xi_M \in K$ for all $M$. Point 2 is clear, but Point 3 requires some work which we now do. 

Since $\xi_M \to \beta^*$ there exists a $\delta > 0$ such that for all $M$,  $\xi_M \in R_{\delta} \subset \mathbb{R}^d$ where $R_{\delta}$ is a square centered at $\beta^*$ with lengths $\delta$. Let $\{v_k\}_{k=1, \cdots, k=2^d}$ be the vertices of $R_{\delta}$. Then since $R_{\delta}$ is a lattice, for any $\xi_M \in R_{\delta}$ we can write

\begin{align}
    \xi_M = \sum_{k=1}^{2^d} \pi_k v_k, \textrm{ where } \sum_{k=1}^{2^d} \pi_i &= 1.
\end{align}

Then by convexity of $s \mapsto e^{s^T x}$ we have
\begin{align}
e^{\xi_M^T(x_i - \bar x)} &= e^{\sum_{k=1}^{2^d} \pi_k v_k^T(x_i - \bar x)} \\
&\leq \sum_{k=1}^{2^d} \pi_k e^{v_k^T (x - \bar x)}.\label{eq:convex_bound}
\end{align}
Let 
\begin{equation}
    d(x) := \sum_{k=1}^{2^d} \pi_k e^{v_k^T (x - \bar x)}(\|x\|^2 + \|\bar x \|^2), 
\end{equation}
which is in $L_{F_0}(\Omega)$ by assumption \eqref{eq:finite_moments}
Then by \eqref{eq:convex_bound} and the triangle inequality we have for all $\xi_M \in K$
\begin{equation} 
e^{\xi_M^T(x - \bar x)}(x^k-\bar x^k)(x^j - \bar x^j) \leq \sum_{k=1}^{2^d} \pi_k e^{v_k^T (x - \bar x)} (\|x \|^2 +\|\bar x \|^2),
\end{equation}
and hence condition 3. is satisfied. We can therefore conclude from the Uniform Law of Large Numbers that $G(N,M)$ converges uniformly in $M$ as $N \to +\infty$. Thus condition 2 of the Moore-Osgood criteria is satisfied and the limit $\lim_{(M,N) \to +\infty}G(M,N)$ exists. Moreover

\begin{align}
\lim_{(M,N) \to +\infty}G(M,N) = \lim_{M \to + \infty} \lim_{N \to +\infty}G(M,N) = \mathbb{E}_{F_0}[e^{\beta_*^T (x-\bar x) }(x-\bar x)^T (x-\bar x)],
\end{align}
almost surely. 
In particular, for $M=N$ we have almost surely
\begin{equation}
    G(N,N) := \frac{1}{N}\sum_{i=1}^N e^{\xi_N^T(x_i - \bar x)}(x_i^k-\bar x^k)(x_i^j - \bar x^j) \to \mathbb{E}_{F_0}[e^{\beta_*^T (x-\bar x) }(x-\bar x)^T (x-\bar x)],
\end{equation}
which yields the result. 
$\Box$

We are now ready to Prove Theorem \ref{theorem:main}.

\noindent \textbf{Proof of Theorem \ref{theorem:main}}

We begin by performing a second order Taylor expansion of \eqref{eq:logloss_empirical2} with respect to $\beta$ around $\beta_*$:
\begin{equation}\label{eq:taylor_exp_orig}
 0 = \nabla_{\beta}\mathcal{ L}_N(\beta_N,\alpha_N) = \nabla_{\beta}\mathcal{L}_N(\beta_*,\alpha_N) +  D_{\beta}^2\mathcal{L}_N(\xi_N,\alpha_N)(\beta_N - \beta_*),
\end{equation}
where $\xi_N = t \beta_N + (1-t) \beta_*$ for some $t \in [0,1]$. 
Since $F_0$ contains $\bar x$ in the sense of Definition 1, we know from Lemma \ref{lemma:alpha_conv} that $e^{\alpha_N} = O(N^{-1})$ and  $\beta_N \to \beta_*$ as $N \to +\infty$ where $\beta_*$ is the unique solution to \eqref{eq:owen_main}. \\

Computing the first two derivatives of \eqref{eq:logloss_empirical2} with respect to $\beta$ evaluated at $\beta_*$ and $\xi_N$ respectively, we have

\begin{align}
    \label{eq:der1a}\nabla_{\beta} \mathcal{L}_N(\alpha_N,\beta_*) &=  -e^{\alpha_N} \sum_{i=1}^N \frac{e^{\beta_*(x_i - \bar x)}}{1+ e^{\alpha_N+\beta_* (x_i-\bar x)}}(x_i-\bar x)\\ &-e^{\alpha_N} \sum_{j=1}^N \frac{e^{\beta_*(x_j - \bar x)}}{1+ e^{\alpha_N+\beta_* (x_j-\bar x)}}(x_j-\bar x) \nonumber \\
    \label{eq:der2a}D_{\beta}^2 \mathcal{L}_N(\alpha_N,\xi_N) &=  -e^{\alpha_N} \sum_{i=1}^N \frac{e^{\xi_N^T(x_i - \bar x)}}{(1+ e^{\alpha_N+\xi_N^T (x_i-\bar x)})^2}(x_i-\bar x) (x_i - \bar x)^T\\ &-e^{\alpha_N} \sum_{j=1}^N \frac{e^{\xi_N^T(x_j - \bar x)}}{(1+ e^{\alpha_N+\xi_N^T (x_j-\bar x)})^2}(x_j-\bar x) (x_j - \bar x)^T \nonumber
\end{align}

By assumption \eqref{eq:xjbound}, the second lines of both \eqref{eq:der1a} and \eqref{eq:der2a} are both $O(e^{\alpha_N})$. Hence \eqref{eq:der1a} and \eqref{eq:der2a} become 
\begin{align}
    \label{eq:der1}\nabla_{\beta} \mathcal{L}_N(\alpha_N,\beta_*) &=  -e^{\alpha_N} \sum_{i=1}^N \frac{e^{\beta_*(x_i - \bar x)}}{1+ e^{\alpha_N+\beta_* (x_i-\bar x)}}(x_i-\bar x) + O(N^{-1})\\
    \label{eq:der2}D_{\beta}^2 \mathcal{L}_N(\alpha_N,\xi_N) &=  -e^{\alpha_N} \sum_{i=1}^N \frac{e^{\xi_N^T(x_i - \bar x)}}{(1+ e^{\alpha_N+\xi_N^T (x_i-\bar x)})^2}(x_i-\bar x) (x_i - \bar x)^T+O(N^{-1})
\end{align}
Since $\beta_N \to \beta_*$, for any $r>0$ we have for sufficiently large $N$,  $\beta_N \in \bar B_{r}(\beta_*)$. Applying Lemma \ref{lemma:sum_approx} to \eqref{eq:der1}-\eqref{eq:der2} along with $e^{\alpha} \leq \frac{cn}{\delta N}$ for some $\delta >0$ from Lemmas \ref{lemma:alpha_conv} and \ref{lemma:der_bounds}, we have a.s

\begin{align}
   \label{eq:first_der_exp} \frac{1}{N}\sum_{i=1}^N \frac{e^{(\beta_*)^T(x_i - \bar x)}}{1+ e^{\alpha_N+(\beta_*)^T (x_i-\bar x)}}(x_i-\bar x) &=  \frac{1}{N}\sum_{i=1}^N e^{(\beta_*)^T(x_i - \bar x)}(x_i-\bar x) + O\left(N^{-1}\right)\\
   \label{eq:second_der_exp} \frac{1}{N}\sum_{i=1}^N \frac{e^{\xi_N^T(x_i - \bar x)}}{(1+ e^{\alpha_N+\xi_N^T (x_i-\bar x)})^2}(x_i-\bar x) (x_i - \bar x)^T &= \mathbb{E}_{F_0}[e^{\beta_*^T (x-\bar x)} (x- \bar x) (x-\bar x)^T]\\ &+ O(N^{-1}),\nonumber
\end{align}
as $N \to +\infty$. \\

Rearranging \eqref{eq:taylor_exp_orig} and substituting in  \eqref{eq:first_der_exp} and\eqref{eq:second_der_exp} we have 
\begin{align}
  \nonumber \sqrt{N} (\beta_N - \beta_*) = - \sqrt{N}\left(D_{\beta}^2 \mathcal{L}_N(\alpha_N, \xi_N )\right)^{-1} &\frac{1}{N}\sum_{i=1}^N\nabla_{\beta} g(x_i; \beta_*)\\ &+ O(\sqrt{N}^{-1}),\label{eq:CLTg}
\end{align}

 Now we claim that $\nabla_{\beta} g(x ; \beta_*)$ has zero expectation. Indeed
\begin{align}
\mathbb{E}_{F_0}[\nabla_{\beta} g(\cdot ; \beta)]= \int e^{\beta (x-\bar x)}(x-\bar x)dF_0(x)=0,
\end{align}
 by \eqref{eq:owen_main}. Since $\nabla_{\beta} g( . ; \beta)$ are i.i.d, by the Central Limit Theorem, we have 
 \begin{equation}\label{eq:clt_num}
     \frac{1}{N}\sum_{i=1}^N\nabla_{\beta} g(x_i; \beta_*) \overset{\mathcal{D}}{\to} \mathcal{N}(0, \textrm{Var}_{F_0}(\nabla_{\beta} g)).
 \end{equation}
 By Lemma \ref{lemma:der_bounds}, we have almost surely that 
 \begin{equation}\label{eq:clt_den}
      D_{\beta}^2 \mathcal{L}_N(\alpha_N, \xi_N) \to \mathbb{E}_{F_0} [ D_{\beta}^2 g(\cdot ; \beta_*)].
 \end{equation}

By Slutskly's Theorem applied to \eqref{eq:CLTg}, using \eqref{eq:clt_num}-\eqref{eq:clt_den}, we conclude that 
\begin{equation}
      \sqrt{N} (\beta_N - \beta_*) \to \mathcal{N}(0, \Sigma_{\bar x, \mu}),
\end{equation}
where 
\begin{equation}
    \Sigma_{\bar x \mu} :=  (\mathbb{E}_{F_0}[\nabla_{\beta}^2 g(\cdot ; t)])^{-1})^T\textrm{Var}_{F_0}(\nabla_{\beta} g(\cdot ; t))(\mathbb{E}_{F_0}[\nabla_{\beta}^2 g(\cdot ; t)])^{-1},
\end{equation}
which proves the Theorem.
\qedsymbol

\section{Simulations}\label{section:simulations}
In order to validate Theorem 1 numerically, we fix $d = 1$, $n=1$, $X \lvert Y=0 \sim \mathcal{N}(0,1)$ and consider $N=100,200,500,1000,5000$. From Owen's result, we know that $\beta_* = \bar x$. Our results are shown in Figure 1.  For each $\bar x$, we obtain the value of the blue curve (empirical) by generating $N$ samples $x_i \sim \mathcal{N}(0,1)$ and solving the Logistic Regression problem using Newtons Method with no regularization. The blue (theorem1) curve is the variance given by Corollary 1, with $\mu=0$ and $\sigma = 1$, namely $\sigma_{\bar x \mu}^2 = e^{\bar x^2}[\bar x^2 + 1]$. \\

We see that our parametric variance from Corollary 1 closely matches the Monte Carlo-based estimate of the variance when $\bar x$ is small in all cases. As $\bar x$ moves further away from $\mu=0$, the approximation deteriorates. But as $N$, the number of samples used for training, increases, the approximation improves and matches closely that of $\sigma_{\bar x, \mu}$.
\begin{figure}[!ht]
  \centering
  {\includegraphics[width=.4\textwidth]{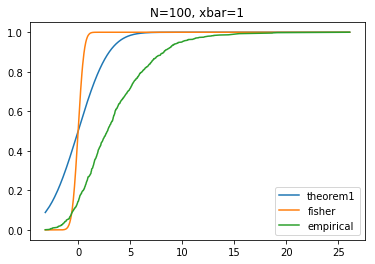}}\quad
  {\includegraphics[width=.4\textwidth]{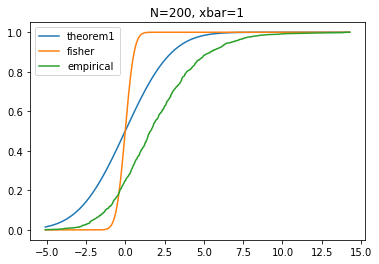}}\\
 {\includegraphics[width=.4\textwidth]{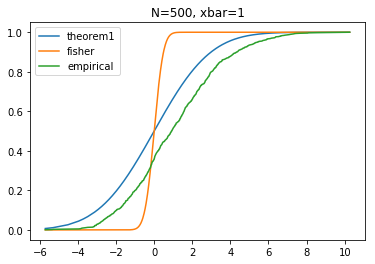}}\quad
{\includegraphics[width=.4\textwidth]{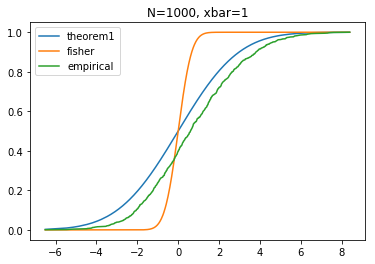}}\\
  {\includegraphics[width=.4\textwidth]{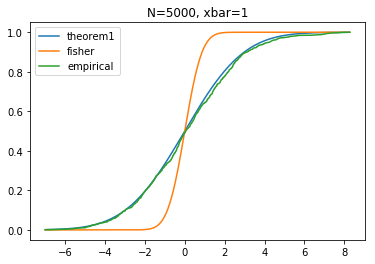}}
  \caption{In each figure we fixed $\bar x = 1$, we have generated $N$ samples from $\mathcal{N}(0,1)$. We solve for for $\beta_N$ using Logistic Regression over 100 Monte Carlo runs and plot the empirical CDF (green) vs the one induced from Theorem 1 (blue).}
  \label{fig:sub1}
\end{figure}

\begin{figure}[!ht]
  \centering
  {\includegraphics[width=.4\textwidth]{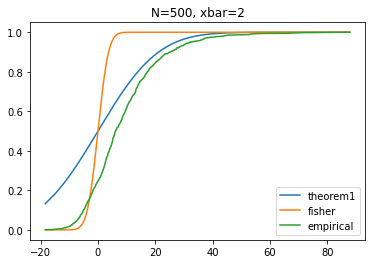}}\quad
  {\includegraphics[width=.4\textwidth]{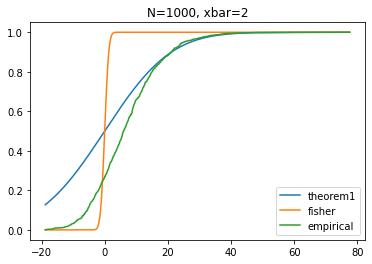}}\\
 {\includegraphics[width=.4\textwidth]{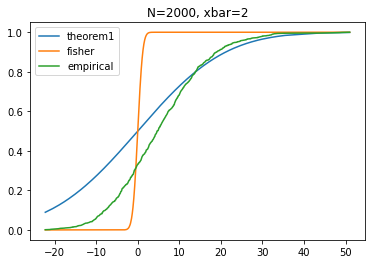}}\quad
{\includegraphics[width=.4\textwidth]{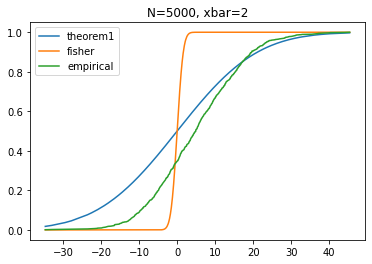}}\\
  {\includegraphics[width=.4\textwidth]{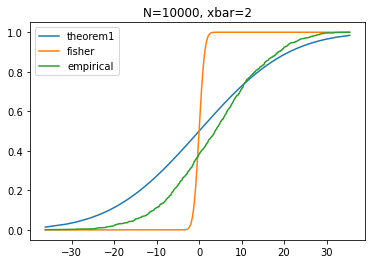}}
  \caption{Here we repeat the analysis in Figure 1 but use $\bar x = 2$. We see that the convergence is significantly slower.}
  \label{fig:sub1}
\end{figure}

\section{Discussion}

In this paper we have derived a second order limit to the limiting slope characterized by  \cite{owen}, and shown that it surprisingly also depends only on the average of the minority class ample. Our result quantifies the extent to which this approximation holds in a way explicitly depending on the distance between $\bar x$ to the center of mass of $F_0$. We have seen that this allows us to fully characterize the asymptotic distribution of the slope parameter in the highly imbalanced case. 
\section{Appendix}

Here we prove some technical lemmas which characterize our result in the case that $F_0$ is normal. 

\begin{lemma}
Let $t \in \mathbb{R}$ and $F_0 = \mathcal{N}(\mu, \sigma^2)$. Then 
\begin{equation}
    \int e^{tx} f(x) dF_0(x) = e^{\frac{t^2\sigma^2}{2}+\mu t } \int f(x) dF_t(x),
\end{equation}
where $F_t = \mathcal{N}(\mu + t\sigma^2, \sigma^2)$.
\end{lemma}
\begin{corollary}\label{cor:variance_ratio}
Let $t = \frac{\bar x- \bar \mu}{\sigma^2}$ and $F_0$ as above. Then 
\begin{align}
    \int e^{tx} f(x) dF_0(x) &= e^{\frac{(\bar x - \bar \mu)^2}{2\sigma^2}+\mu\frac{ (\bar x - \bar \mu)}{\sigma^2} } \int f(x) dF_t(x)\\
    \int e^{2tx} f(x) dF_0(x) &=  e^{\frac{4(\bar x - \bar \mu)^2}{2\sigma^2}+2\mu \frac{(\bar x - \bar \mu)}{\sigma^2} } \int f(x) dF_{2t}(x).
\end{align}
In particular, 
\begin{equation}
    \frac{ \int e^{2tx} f(x)dF_0(x)}{\left(\int e^{tx} f(x) dF_0(x)\right)^2} =  e^{\frac{(\bar x - \bar u)^2}{\sigma^2}} \frac{\int f(x) dF_{2t}(x)}{\left(\int f(x) dF_t(x)\right)^2}
\end{equation}
\end{corollary}

We now focus on deriving an explicit expression for $\sigma_{\bar x, \mu}^2$. Using Corollary \ref{cor:variance_ratio} we conclude that 
\begin{equation}
    \sigma_{\bar x \mu}^2 = e^{\frac{(\bar x - \mu)^2}{\sigma^2}} \frac{\int (x - \bar x)^2 dF_{2(\bar x - \mu)/\sigma^2}(x)}{\left(\int (x - \bar x)^2 dF_{(\bar x - \mu)/\sigma^2}(x)\right)^2}
\end{equation}
 Evaluating explicitly, we obtain 
\begin{equation}
     \sigma_{\bar x \mu}^2 = e^{\frac{(\bar x - \mu)^2}{\sigma^2}}\left[\frac{(\bar x - \mu)^2 + \sigma^2}{\sigma^4}\right].
\end{equation}
This characterizes how errors in the estimates depend crucially on the distance between the support of the points where $y=1$ and $y=0$ - it is in fact exponential. This gives us an explicit calculation we can perform to determine if the Owen approximation will suffice given the size of our training data.

These final Lemmas are used in the proof of Theorem 1.
\begin{lemma}
Let $A(n) = [a_{ij}^n]_{i,j=1,\cdots,M}$. Then if $A(n) \to A$ and $A$ is invertible, then $A(n)$ is invertible for sufficiencly large $n$. and $A(n)^{-1} \to A^{-1}$ 
\end{lemma}
\noindent \textbf{Proof:} This is an immediate consequence of the continuity of the determinant and the fact that $A$ is invertible iff $\det(A)>0$. $\qedsymbol$\\

\begin{lemma}
Let $Z \sim \mathcal{N}(\mu, \Sigma)$ with $\mu \in \mathbb{R}^d$ and $\Sigma \in \mathbb{R}^d \times \mathbb{R}^d$ Then $AZ \sim \mathcal{N}(A \mu, A \Sigma A^T)$.
\end{lemma}
\noindent \textbf{Proof:} If $Z$ has mean $\mu$, then $X = AZ$ has mean $\mathbb{E}[X] = A\mu$ and covariance $\mathbb{E}[(X-A\mu)^T (X-A \mu)] = \mathbb{E}[(Az - A\mu)^T (Az-A\mu)] = \mathbb{E}[(z-\mu)^TA^TA  (z-\mu)] =  \mathbb{E}[A(z-\mu)  (z-\mu)^T A^T] =A \Sigma A^T$.
$\qedsymbol$\\

\medskip
{
\small
\bibliographystyle{plainnat}
\bibliography{lr}

\begin{thebibliography}{8}
\providecommand{\natexlab}[1]{#1}
\providecommand{\url}[1]{\texttt{#1}}
\expandafter\ifx\csname urlstyle\endcsname\relax
  \providecommand{\doi}[1]{doi: #1}\else
  \providecommand{\doi}{doi: \begingroup \urlstyle{rm}\Url}\fi

\bibitem[Agresti(2015)]{agresti2015foundations}
Alan Agresti.
\newblock \emph{Foundations of linear and generalized linear models}.
\newblock John Wiley \& Sons, 2015.

\bibitem[Billingsley(1999)]{billingsley1999convergence}
Patrick Billingsley.
\newblock \emph{Convergence of Probability Measures}.
\newblock John Wiley \& Sons, 1999.

\bibitem[Bolton et~al.(2002)Bolton, Hand, et~al.]{bolton2002statistical}
Richard~J Bolton, David~J Hand, et~al.
\newblock Statistical fraud detection: A review.
\newblock \emph{Statistical science}, 17\penalty0 (3):\penalty0 235--255, 2002.

\bibitem[Fan(2016)]{ZF2016}
Zhou Fan.
\newblock Statistics 200: Introduction to statistical inference, 2016.
\newblock URL
  \url{https://web.stanford.edu/class/archive/stats/stats200/stats200.1172/Lecture16.pdf}.

\bibitem[Lee et~al.(2012)Lee, Orten, Dasdan, and Li]{lee2012estimating}
Kuang-chih Lee, Burkay Orten, Ali Dasdan, and Wentong Li.
\newblock Estimating conversion rate in display advertising from past
  erformance data.
\newblock In \emph{Proceedings of the 18th ACM SIGKDD international conference
  on Knowledge discovery and data mining}, pages 768--776, 2012.

\bibitem[Owen(2007)]{owen}
Art~B. Owen.
\newblock Infinitely imbalanced logistic regression.
\newblock \emph{Journal of Machine Learning Research}, 8:\penalty0 761--773,
  2007.

\bibitem[Silvapulle(1981)]{Silvapulle1981}
Mervyn~J. Silvapulle.
\newblock On the existence of maximum likelihood estimators for the binomial
  response models.
\newblock \emph{Journal of the Royal Statistical Society. Series B
  (Methodological)}, 43\penalty0 (3):\penalty0 310--313, 1981.
\newblock ISSN 00359246.
\newblock URL \url{http://www.jstor.org/stable/2984941}.

\bibitem[Zhu et~al.(2006)Zhu, Su, and Chipman]{zhu2006lago}
Mu~Zhu, Wanhua Su, and Hugh~A Chipman.
\newblock Lago: A computationally efficient approach for statistical detection.
\newblock \emph{Technometrics}, 48\penalty0 (2):\penalty0 193--205, 2006.

\end{thebibliography}
}




\end{document}